\documentclass[article]{siamart171218}
\usepackage{amsmath}
\usepackage{float}
\usepackage{diagbox}
\usepackage{bbold}
\usepackage{hyperref}
\usepackage{float}
\usepackage{diagbox}
\usepackage{bbold}
\usepackage{graphicx}
\usepackage{wrapfig}

\newcommand{\ep}{\epsilon} 
\newcommand{\Aep}{A^{\epsilon} }
\newcommand{\bep}{b^{\epsilon} }
\newcommand{\cep}{c^{\epsilon}}
\newcommand{\Fep}{F^{\epsilon} }
\newcommand{\dt}{\Delta t }
\newcommand{\At}{\tilde{A} }
\newcommand{\bt}{\tilde{b} }
\newcommand{\ct}{\tilde{c} }

\begin{document}

\title{Perturbed Runge-Kutta methods for mixed precision applications\thanks{Submitted to the editors on December 24, 2020
\funding{This material is based upon work supported by the U.S. Department of Energy, Office of Science, Office of Advanced Scientific Computing Research, as part of their Applied Mathematics Research Program. The work was performed at the Oak Ridge National Laboratory, which is managed by UT-Battelle, LLC under Contract No. De- AC05-00OR22725. The United States Government retains and the publisher, by ac- cepting the article for publication, acknowledges that the United States Government re- tains a non-exclusive, paid-up, irrevocable, world-wide license to publish or reproduce the published form of this manuscript, or allow others to do so, for the United States Government purposes. The Department of Energy will provide public access to these results of federally sponsored research in accordance with the DOE Public Access Plan 
\url{http://energy.gov/downloads/doe-public-access-plan}}}}

\author{Zachary J. Grant \thanks{Department of Computational and Applied Mathematics, Oak Ridge National Laboratory, Oak Ridge TN 37830.       \email{grantzj@ornl.gov}}}

\maketitle

\begin{abstract}
In this work we consider a mixed precision approach to accelerate the implemetation of multi-stage methods.
We show that Runge--Kutta methods can be designed so that certain costly intermediate computations
can be performed as a  lower-precision computation without adversely impacting the accuracy of 
the overall solution. In particular, a properly designed Runge--Kutta method  will damp out the 
errors committed in the initial stages. This is of particular interest when we consider implicit Runge--Kutta
methods. In such cases, the implicit computation of the stage values can be considerably faster if the
solution can be of lower precision (or, equivalently, have a lower tolerance). We provide a general theoretical
additive framework  for designing mixed precision  Runge--Kutta methods, and use this framework
to derive order conditions for such methods. Next, we show how using this approach allows us
to leverage low precision computation of the implicit solver while retaining high precision in the overall method.
We  present the behavior of some mixed-precision implicit Runge--Kutta methods through numerical studies,
and demonstrate how the numerical results match with the theoretical framework. This novel 
mixed-precision implicit Runge--Kutta framework opens the door to the design of many such methods.
\end{abstract}

\section{Introduction}
Consider the ordinary differential equation (ODE)
\[ u_t = F(u) .\]
Evolving this equation using a standard Runge Kutta approach we have the $s$-stage Runge--Kutta method
\begin{eqnarray}
y^{(i)}= u^n + \dt \sum_{j=1}^{s} A_{ij} F(y^{(j)}) \nonumber \\
u^{n+1}=u^n + \dt \sum_{j=1}^{s} b_j F(y^{(j)}),
\end{eqnarray}
where  $A^{s x s}$ and  $b^{1 x s} $  are known as the Butcher coefficients of the method.

The computational cost of the function evaluations $F(y^{(j)}) $ can be considerable, especially in cases where 
it necessitates an implicit solve. The use of mixed precision approach has been implemented for other  numerical methods
 \cite{Abdelfattah2020, WENOMP} seems to be a promising approach. 
Lowering the precision on these computations, either by storing $F(y^{(j)})$ as a
single precision variable rather than a double precision one, or by raising the tolerance of the implicit solver,
can speed up the computation significantly. However, it generally reduces the precision of the overall numerical solution.

The aim of this paper is to create a framework for the design of Runge--Kutta methods that allow lower precision
function evaluations for some of the stages, without impacting the overall precision of the solution. 
This acceleration of multi-stage methods using a mixed precision approach relies on a design that ensures 
that errors committed in the early stages may be damped out by the construction of the update in later stages.

The structure of this paper is as follows: in Section \ref{sec:motivation} we provide a numerical example
of a mixed-precision formulation of the implicit midpoint rule, which motivates the need to study the effect
of lower-precision computations of implicit function evaluations.
In  \ref{sec:framework} we provide a general framework for analyzing mixed-precision Runge--Kutta methods 
by exploiting the additive Runge--Kutta method formulation. In this section we also present the order conditions
that arise from this formulation and show how these are a relaxation of the general required order conditions.
In Section \ref{sec:methods} we show how the implicit midpoint rule  can be written in this additive 
mixed-precision Runge--Kutta  form, and develop methods in a specific class of implicit Runge--Kutta methods that 
exploit the formulation in Section \ref{sec:framework}. 

\section{Motivating example: mixed precision implementation of the implicit midpoint rule} \label{sec:motivation}

First we consider the implicit midpoint rule defined as
\begin{align}\label{IMPClassic}
u^{n+1} = u^n + \Delta t F\left(\frac{u^n+u^{n+1} }{2} \right)
\end{align}
which is equivalently written in its Butcher form as:
\begin{subequations} \label{IMPButcher}
\begin{align}
y^{(1)}&= u^n + \frac{\Delta t}{2} F\left(y^{(1)} \right) \label{IMPButcher1}\\
u^{n+1}&= u^n + \Delta t F \left(y^{(1)} \right) \label{IMPButcher2}
\end{align}
\end{subequations}
The equivalence can be seen by noting that equation \eqref{IMPButcher2} can be re-written as $u^{n+1}= 2y^{(1)}-u^n$, 
and substituting $y^{(1)} = \frac{1}{2} (u^n + u^{n+1})$  into equation \eqref{IMPButcher1}.

  In both equations \eqref{IMPClassic} and \eqref{IMPButcher1}, we require an implicit solver  in order to compute the update.  
  Let's consider the case that the implicit solver can only be satisfied up to some tolerance, $O(\ep)$, who's output ${u_\ep}^{n+1}$ satisfies the modified equation
  \[ {u_\ep}^{n+1} = u^n + \Delta t \Fep \left( \frac{u^n+{u_\ep}^{n+1} }{2} \right) \]
  exactly. We define the perturbation operator 
  \begin{eqnarray} \label{tau}
   \tau(u) =\frac{F(u)-\Fep(u)}{\ep}  
   \end{eqnarray}
  so at any given time-step  $u^{n+1}- {u_\ep}^{n+1}$ = $O( \ep \dt)$. 
  Over the course of the solution, the local errors build-up to give an global error contribution of  $O(\ep)$.
  
  However, we can also formulate a method
  \begin{subequations} \label{IMPButcherMP}
  \begin{align}
{y_\ep}^{(1)}&= u^n + \frac{\Delta t}{2} \Fep \left(y_\ep^{(1)} \right)  \label{IMPButcherMP1}\\
\hat{u}^{n+1}&= u^n + \Delta t F\left( y_\ep^{(1)}\right) \label{IMPButcherMP2}
\end{align}
\end{subequations}
  where we use an inaccurate implicit solve in the first stage, and an accurate explicit evaluation in the second stage.
We obtain $y_\ep^{(1)} = y^{(1)} + \dt O(\ep)$.   The error, ${u}^{n+1} - \hat{u}^{n+1}$,  between the full precision and mixed precision implementation over
 one time-step  is now of order   $O(\ep \dt^2 )$. 
  This results in a  global error contribution of $O(\ep \dt )$.

  \begin{wrapfigure}[24]{r}{0.55\linewidth} \vspace{-0.05in}
\includegraphics[width=0.55\textwidth]{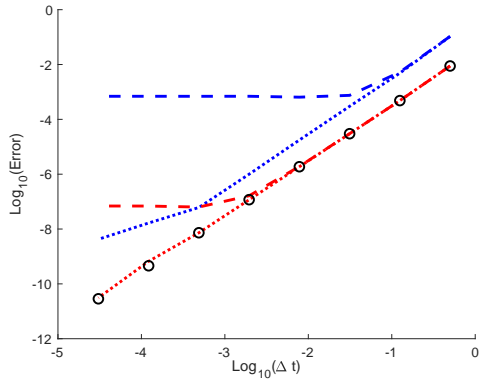}
\caption{{The implicit midpoint rule applied to the van der Pol system \eqref{eqn:vdp}.
The half-precision implementation ($\ep= O(10^{-4})$) is shown as a blue dashed line 
and the single-precision implementation  $\ep= O(10^{-8})$ as a red dashed line. 
The mixed precision implementation with $\ep= O(10^{-4})$ is shown with a dotted blue line, 
and the mixed precision implementation with $\ep= O(10^{-8})$ is shown with a dotted red line. 
The double precision implementation is shown for reference with  black circles.
}}
\label{fig:IMR1}
\end{wrapfigure} To see how this works in practice, consider the van der Pol system
  \begin{subequations} \label{eqn:vdp} \vspace{-.2in}
  \begin{eqnarray}
  y_1' & = &  y_2 \\
  y_2' & = &  \; y_2  \; (1-y_1^2) - y_1
  \end{eqnarray}
\end{subequations}
  with initial conditions $y_1(0)= 2$ and $y_2(0) = 0$. We stepped this forward to a final time $T_f=1.0$. 
  To emulate a  mixed precision computation,  we let  $\Fep$  in  \eqref{IMPButcherMP1} be the truncated output of $F$, 
  where the truncation is performed  to single precision  ($\ep \approx 10^{-8}$) and half precision  ($\ep \approx 10^{-4}$),
  while the explicit evaluation of $F$ in   \eqref{IMPButcherMP2} is performed in double precision.
For comparison, we emulate  a low precision  version of \eqref{IMPButcher}, where we truncate all computations of 
$F$ in to the same  low precision.  Finally, we compute a full -precision version  of \eqref{IMPButcher}, 
where all $F$ are computed to double precision and not truncated.
(Note that his truncation approach is advocated in \cite{Higham2020} to emulate low precision simulations.)

In Figure  \ref{fig:IMR1}, we show the final time errors in the simulation using the various implementations.
First, we look at the errors from a low-precision implementation \eqref{IMPButcher} in which  
{\em all}  the computations of $F$ are truncated to single precision (blue dashed line) or half precision (red dashed line).
Clearly, for a sufficiently refined $\dt$, these errors look like $O(\ep)$. Next, we look at the errors from the mixed precision computation,  
given by Equation \eqref{IMPButcherMP},  has a value of $\Fep $ that results from a single precision (blue dotted line) 
truncation  or a half-precision truncation (red dotted line), and  an explicit $F$ that is evaluated in full precision. 
These errors are clearly much better: they are convergent initially at second order, and the implementation 
with $\ep \approx 10^{-8}$ remains second order,  and matches the full-precision implementation (black circles).
However, the implementation using $\ep \approx 10^{-4}$ reduces to first order for a sufficiently  small $\dt$.
  In the next section we construct a general framework that explains why this happens and allows us to construct higher order
  methods that work well in a mixed precision formulation.

\section{A general framework for the analysis of mixed precision Runge--Kutta methods}  \label{sec:framework} 

We use the B-series analysis  for additive Runge Kutta methods to develop consistency and perturbation  conditions for the 
mixed precision Runge Kutta method of the form:
\begin{subequations}
\begin{align}
y^{(i)}= u^n + \dt \sum_{j=1}^{s}A_{ij}F(y^{(j)}) +  \dt \sum_{j=1}^{s}\Aep_{ij}\Fep (y^{(j)})\\
u^{n+1}= u^n + \dt \sum_{j=1}^{s}b_{j}F(y^{(j)}) +  \dt \sum_{j=1}^{s}\bep_{j}\Fep (y^{(j)}).
\end{align}
\end{subequations}
Such methods have been extensively studied, including \cite{Ascher1995,Carpenter2003,Sandu2015}

In  \cite{KetchPerturbation}  a perturbation approach to  Runge--Kutta methods was proposed  in order to  increase the largest allowable  time-step that preserved strong stability.
Following a similar approach we note that the operator $\Fep(y)$  is an approximation to  $F(y)$ such that  for any $y$ we require  $\Fep(y)-F(y)=O(\ep)$. 
This allows us to rewrite the scheme to evolve the operator $F$ and its perturbation $\tau$, defined in \eqref{tau}:
\begin{subequations} 
\begin{align}
\label{Perturbed Methoda}
y^{(i)}&= u^n + \dt \sum_{j=1}^{s}\tilde{A}_{ij}F(y^{(j)}) +  \ep \dt \sum_{j=1}^{s}\Aep_{ij} \tau(y^{(j)})\\
u^{n+1}&= u^n + \dt \sum_{j=1}^{s}\tilde{b}_{j}F(y^{(j)}) + \ep  \dt \sum_{j=1}^{s} \bep_{j} \tau (y^{(j)}) \label{Perturbed Methodb}
\end{align}
\end{subequations}
where $\tilde{A}_{ij}= A_{ij}+\Aep_{ij}$ and $\tilde{b} = b_{j}+\bep_{j}$.

Analyzing the scheme in this form allows us to use an additive B-series representation to track the evolution of $ F$ 
as well as its interaction with  the perturbation  function $\tau$. For example, a second order expansion is:
\begin{align} \label{expansion}
u^{n+1} &= \underbrace{ u^n +  \dt\bt eF(u^n) + \dt^2\bt \ct F_y(u^n)F(u^n)}_{scheme} \\
& +   \underbrace{\ep\dt \left( \bep e\tau(u^n) +    \dt \left( \bep \ct \tau_y(u^n)F(u^n)  +\bt \cep F_y(u^n) \tau(u^n)
 +  \ep  \bep \cep \tau_y(u^n)\tau(u^n) \right) \right) }_{perturbation}
\nonumber  \\
& + O(\dt^3). \nonumber 
\end{align} 
This expansion shows two sources of error: those of the scheme and those of the perturbation, thus the 
error at one time-step can be written as the sum of the approximation error of the scheme,   $E_{sch} $,
and  the perturbation error, $E_{per}.$
This leads to two sets of conditions under which  a perturbed scheme has an error
$$E= E_{sch} + E_{per} =  O(\dt^{p+1}) + O(\ep \dt^{m+1} )$$ 
at each time-step. A method with these errors will have a global error of the form
\[ Error = O(\dt^{p}) + O( \ep \dt^{m}).\]

We can easily extend \eqref{expansion} to higher order. The fourth order expansion related to the consistency of the scheme has following terms:
\[ \begin{array}{lll} \hline
\mbox{Terms involving  $\dt$} \; \; \; \;  \; \; \; \; & \mbox{Terms involving $F$ } \; \; \; \; \; \; \; & \mbox{scheme} \\ 
\; \; \; \;  & \ & \mbox{coefficients} \\ \hline
& & \\
\dt & F(u^n)  & \bt e \\
\dt^2  & F_y(u^n) F(u^n) & \bt \ct \\ 
\dt^3  & F_y(u^n)F_y(u^n)F(u^n) & \bt \At \ct   \\
\dt^3  &  F_{yy}(u^n)(F(u^n),F(u^n))& \bt (\ct \cdot \ct)  \\
\dt^4  &F_y(u^n)F_y(u^n)F_y(u^n)F(u^n) &  \bt\At \At \ct  \\
\dt^4  &F_y(u^n)F_{yy}(u^n)(F(u^n),F(u^n)) & \bt\At(\ct \cdot \ct) \\
\dt^4  &F_{yy}(u^n)(F_y(u^n)F(u^n),F(u^n)) & \bt(\At \ct \cdot \ct)  \\
\dt^4  &F_{yyy}(u^n)(F(u^n),F(u^n),F(u^n)) & \bt(\ct \cdot \ct \cdot \ct) .
\end{array} \]

Expanding the terms related to the perturbation error to third order in $\dt$ and third order in $\epsilon$ we obtain:
\[ \begin{array}{lll} \hline 
\mbox{Terms involving } \; \; \; \;  & \mbox{Terms involving } & \mbox{scheme} \\ 
\mbox{ $\ep$ and $\dt$} \; \; \; \;  & \mbox{ $F$ and $\tau$} & \mbox{coefficients} \\ \hline
& & \\
\ep \dt \; \; \; \; &  \tau(u^n)      & \bep e \\ 
\ep \dt^2 & \tau_y(u^n)F(u^n)    & \bep \ct \\
\ep \dt^2 & F_y(u^n)\tau(u^n)    & \bt \cep \\ 
\ep^2\dt^2  & \tau_y(u^n)\tau(u^n)& \bep \cep \\
\ep \dt^3 & \tau_y(u^n)F_y(u^n)\tau(u^n)    & \bep \At\ct \\ 
\ep \dt^3 & F_y(u^n)\tau_y(u^n)F(u^n)   & \bt \Aep \ct   \\ 
\ep \dt^3 & F_y(u^n)F_y(u^n)\tau(u^n)    &  \bt \At \cep  \\
\ep \dt^3 &  \tau_{yy}(u^n)(F(u^n),F(u^n))   &  \bep(\ct \cdot \ct) \\
\ep \dt^3  & F_{yy}(u^n)(F(u^n),\tau(u^n))    &  \bt(\ct \cdot \cep)  \\
\ep^2 \dt^3  & \tau_y(u^n)\tau_y(u^n)F(u^n)  \; \; \; \; \; \;& \bep \Aep \ct  \\
 \ep^2 \dt^3  & \tau_y(u^n)F(u^n)\tau(u^n)      & \bep \At \cep \\
  \ep^2 \dt^3  & F_y(u^n)\tau_y(u^n)\tau(u^n)      &  \bt \Aep \cep   \\
    \ep^2 \dt^3  & \tau_{yy}(u^n)(\tau(u^n),F(u^n))    &    \bep(\cep \cdot \ct)  \\
     \ep^3\dt^3 & \tau_y(u^n)\tau_y(u^n)\tau(u^n)      & \bep \Aep \cep \\
     \ep^3 \dt^3 \; \; \; \; \;  \; \; \; \; \; \; &  \tau_{yy}(u^n)(\tau(u^n),\tau(u^n))  \;  \; \; \;    \; \; \; \; \;  \; \; \;   & \bep(\cep \cdot \cep)  \\ \hline 
\end{array} \]
From these terms, the  consistency conditions for the scheme and the perturbation conditions can be easily defined, and are given in the next section.
We consider two possible scenarios: In the standard scenario, we assume that both $F$ and $\Fep$ are well-behaved functions,
 and all of their derivatives exist and are bounded. In this case, for the perturbation
 terms in the table above to be zeroed out, we simply impose the condition that the corresponding coefficient are zero.
 In the mixed precision scenario which motivated this formulation, we  consider that  $\tau$ comes from a precision error 
 that is defined by "chopping" the values at the  desired precision. In this case, $\Fep=Chop(F)$,  so that the 
operator $\Fep$  is bounded but not Lipshitz continuous, and so $\tau_y = \frac{F_y - \Fep_y}{\epsilon}$  is also not Lipshitz continuous,
{\em i.e} $\frac{ \partial^{(k)} \tau}{\partial y^{(k)}} $ does not exist.
 In this case, we must ensure that all terms containing $\tau_y$ are multiplied by zero coefficients in the expansion,
 without assuming cancellation errors.
In other words,  we  requires more stringent conditions to ensure  that terms of the form $\tau_y , \tau_{yy}$  etc. do
not appear in the final stage.  This means that whereas when $\tau$ is a well-behaved function, it is sufficient to require that
$  \bep c=0$, instead we must require that not only the sum is zero, but every term: $ \bep_j c_j=0$.
We denote conditions of this form with absolute values (e.g. $ |\bep| |c| = 0$). These stringent conditions apply to each
coefficient matrix/vector which appears in conditions corresponding to derivatives of tau.
These conditions are presented in the next subsection.

\subsection{Consistency Conditions and Perturbation Conditions} \label{sec:conditions}
The consistency  and perturbation conditions can be derived from the additive B-series analysis of the 
scheme defined by equations (\ref{Perturbed Methoda}-\ref{Perturbed Methodb}). 
Terms that involve only $F$ yield the classical consistency conditions of the unperturbed scheme. 
The cross terms are the perturbation errors.
A perturbed Runge-Kutta method is of consistency order $p$ if the following conditions are satisfied:
\begin{subequations} \label{OC}
\[ \begin{array}{ll}
 \text{For } p\geq1: \; \; \; \;  &   \; \; \; \;   \bt e=1 \\
 \text{For } p\geq2: &  \; \; \; \;  \bt\ct=\frac{1}{2} \\
 \text{For } p\geq 3: & \; \; \; \;   \bt(\ct \cdot \ct) = \frac{1}{3} \; \; \; \mbox{and} \; \; \; \bt\At\ct = \frac{1}{6} \\
\text{For } p\geq 4: & \; \; \; \;  \bt(\ct \cdot \ct \cdot \ct) = \frac{1}{4} , \; \;
 \bt(\At \ct \cdot \ct) = \frac{1}{8}, \; \; \; \; 
 \bt\At(\ct \cdot \ct) = \frac{1}{12}, \; \; \mbox{and} \; \; \;  \bt\At \At \ct = \frac{1}{24} .
 \end{array} \]
 \end{subequations}

The perturbation errors are determined by both $\dt$ and $\ep$. 
For a scheme to achieve order $O( \ep \dt^{m+1})$ we require :
\begin{subequations} \label{PC}
\begin{itemize}
\item For $m\geq 1$ we require 
 \begin{align}
\mbox{for} \; \;  n\geq 1:  \; \; \; \; &  \bep e=0 
 \end{align}
\item For $m\geq 2$ we require 
 \begin{align}
\mbox{for} \; \;  n\geq 1: \; \; \; & |\bep| |\ct| =0 \; \; \mbox{and} \; \;  \bt \cep = 0 \\
\mbox{for} \; \;  n\geq 2:  \; \; \; &  |\bep| |\cep|=0 
 \end{align}

 \item For $m\geq 3$ we require 
 \begin{align}
\mbox{for} \; \;  n\geq 1: \; \; \; \;  & |\bep| |\At| |\ct|=0 , \; \;
| \bt| | \Aep| |\ct|  = 0,  \; \; 
 \bt \At \cep=0, \\
 & 
| \bep| \left|(\ct \cdot \ct)\right|=0, \; \; 
 \; \;  \bt(\ct \cdot \cep) =0, \nonumber  \\
\mbox{for} \; \;  n\geq 2: \; \;  \; \; &  
| \bep| | \Aep|  |\ct|=0, \; \;
| \bep| |\At| |\cep|=0, \; \;
 | \bt| | \Aep| | \cep| = 0, \\ 
& | \bep| \left| (\cep \cdot \ct)\right| =0, \; \;
 \bt(\cep \cdot \cep) =0, \nonumber \\
\mbox{for} \; \;  n\geq 3: \; \; \; \;  & 
|\bep| | \Aep| | \cep| =0 , \; \;   | \bep| \left| (\cep \cdot \cep)\right| =0.
 \end{align}
 \end{itemize}
 \end{subequations}
    Note that  the coefficients that are  not attached to a derivative (see table above),  do not require the absolute value.

In the Subsection \ref{IMRexplained}  we will show how these conditions can explain the behavior of the mixed-precision implementation of the implicit midpoint rule.
In Section \ref{sec:methods} we will use these conditions to derive mixed-precision methods.

\subsubsection{Perturbation conditions when $\tau$ is well behaved}
If $\tau$ is a well-behaved function, we can assume that terms with similar terms will cancel. In this case, the perturbation conditions 
simplify.  For a scheme to achieve order $O(\ep \dt^{m+1} )$ we require:
\begin{subequations} \label{PCsimplified}
\begin{itemize}
\item For $m\geq 1$ we require 
 \begin{align}
\mbox{for} \; \;  n\geq 1:  \; \; \; \; &  \bep e=0 
 \end{align}
\item For $m\geq 2$ we require 
 \begin{align}
\mbox{for} \; \;  n\geq 1: \; \; \; & \bep \ct=0 \; \; \mbox{and} \; \;  \bt\cep = 0 \\
\mbox{for} \; \;  n\geq 2:  \; \; \; &  \bep \cep=0 
 \end{align}
 \item For $m\geq 3$ we require 
 \begin{align}
\mbox{for} \; \;  n\geq 1: \; \; \; \;  & \bep \At\ct=0 , \; \;
 \bt\Aep\ct = 0,  \; \; 
 \bt\At\cep=0, \\
 & 
 \bep(\ct \cdot \ct)=0, \; \; 
 \; \;  \bt(\ct \cdot \cep) =0, \nonumber  \\
\mbox{for} \; \;  n\geq 2: \; \;  \; \; &  
\bep \Aep \ct=0, \; \;
 \bep \At \cep=0, \; \;
  \bt \Aep \cep = 0, \\ 
&  \bep(\cep \cdot \ct)=0, \; \;
 \bt(\cep \cdot \cep) =0, \nonumber \\
\mbox{for} \; \;  n\geq 3: \; \; \; \;  & 
\bep \Aep \cep=0 , \; \;    \bep(\cep \cdot \cep) =0.
 \end{align}
 \end{itemize}
 \end{subequations}

\subsection{Understanding the mixed precision implicit midpoint rule using the additive framework} \label{IMRexplained}
Returning to the implicit midpoint rule example in Section \ref{sec:motivation}, 
we can use the framework developed above  to study  the numerical behavior of the different implementations
of the implicit midpoint rule. First, we consider the high precision form \eqref{IMPButcher}:
\begin{subequations}
\begin{align*}
y^{(1)}&= u^n + \frac{\Delta t}{2} F\left(y^{(1)} \right)\\
u^{n+1}&= u^n + \Delta t F \left(y^{(1)} \right) 
\end{align*}
\end{subequations}
and note that it is exactly equivalent to \eqref{IMPClassic}. 
The low precision form of  \eqref{IMPButcher} and \eqref{IMPClassic} is 
\begin{subequations} \label{IMPButcher_LP}
\begin{align} 
y^{(1)}&= u^n + \frac{\Delta t}{2} {\Fep}\left(y^{(1)} \right)\\
u^{n+1}&= u^n + \Delta t \Fep\left(y^{(1)} \right) 
\end{align}
\end{subequations}
while the mixed precision form is
\begin{subequations} \label{IMPButcher_MP}
\begin{align} 
y^{(1)}&= u^n + \frac{\Delta t}{2} \Fep\left(y^{(1)} \right)\\
u^{n+1}&= u^n + \Delta t F \left(y^{(1)} \right) .
\end{align}
\end{subequations}
Now let's look at the coefficients of each of these schemes.
The full-precision method  \eqref{IMPButcher}  has coefficients:
\[    
\hspace{.2in}
b = \left( \begin{array}{l}
1\\
0 \\
\end{array} \right) ,
\hspace{.2in}
c = \left( \begin{array}{l}
 \frac{1}{2} \\
1 \\
\end{array} \right) 
\hspace{.2in}
\bep = \left( \begin{array}{l}
0\\
0 \\
\end{array} \right) ,
\hspace{.2in}
\cep = \left( \begin{array}{l}
0 \\
0 \\
\end{array} \right) 
 \]
 which satisfy the order conditions
\[ \bt e = be+ \bep e  = 1, \; \; \; \;
 \bt  \ct = bc  + b\cep + \bep c +  \bep \cep = \frac{1}{2}.\]
The errors from the reduced precision operator are all zero: 
\[ 
\bep e  = 0, \; \; \; \; 
 b\cep = \bep c = \bep\cep  = 0, \]
so that, as  expected,  there is no low precision contribution. This 
 high precision method will produce second order global errors:  $ Error = O(\dt^2)$.

Next, we look at the low-precision method \eqref{IMPButcher_LP}
\[  
b = \left( \begin{array}{l}
0\\
0 \\
\end{array} \right), 
\hspace{.2in}
c = \left( \begin{array}{l}
0\\
0\\
\end{array} \right),
\hspace{.2in}
\bep = \left( \begin{array}{l}
1\\
0 \\
\end{array} \right),
\hspace{.2in}
\cep = \left( \begin{array}{l}
\frac{1}{2} \\
1  \\
\end{array} \right).
 \]
Once again the $O(\dt^2)$  consistency  conditions are satisfied
\[\bt e = be+ \bep e  = 1, \; \; \; \;
\bt \ct = bc  + b\cep + \bep c +  \bep \cep = \frac{1}{2} .\]
However,  the errors introduced by the reduced precision operator are given by the perturbation condition
\begin{eqnarray*}
\bep e  =  1  ,
\end{eqnarray*}
so that at each time-step we have  a perturbation error of the form
\[ E_{per} = \ep\dt\bep e \tau(u^n) = \dt O(\ep). \]
Putting this together, we expect to see errors of the form 
\[ E = O(\dt^3) + O(\dt  \epsilon) \]
at each time-step, and an overall error of 
\[Error =O(\dt^2) + O( \epsilon) . \]
This explains the reduced accuracy we see  in Figure \ref{fig:IMR1},
where initially the errors of $O(\dt^2)$ dominate, but as  $\dt$ gets small enough, the $O( \epsilon) $
terms dominate.

Finally, we look at the mixed-precision method \eqref{IMPButcher_MP}
\[  
b = \left( \begin{array}{l}
1\\
0 \\
\end{array} \right),
\hspace{.2in}
c = \left( \begin{array}{l}
0\\
1\\
\end{array} \right) ,
\hspace{.2in}
\bep = \left( \begin{array}{l}
0\\
0 \\
\end{array} \right) ,
\hspace{.2in}
\cep = \left( \begin{array}{l}
\frac{1}{2} \\
0  \\
\end{array} \right) .
 \]
 We observe that, as before, the $O(\dt^2)$ consistency conditions are satisfied
\[  \bt e = be+ \bep e  = 1,   \; \; \; \;
\bt \ct = bc  + b\cep + \bep c +  \bep \cep = 0 + \frac{1}{2} + 0 + 0 = \frac{1}{2} ,
\]
so that $E_{sch} = \dt^3$.
The perturbation  errors from the reduced precision operator are
\[  \bep e  =  0  , \; \; \; \;
 b \cep = \frac{1}{2},  \; \; \bep c=0,  \; \;   \bep\cep = 0 , \]
 so that at each step we will see perturbation errors of the form \[ E_{per} = \ep \; \dt^2 b \cep \; F'(u^n) \tau(u^n) =\frac{1}{2} \;  \ep \dt^2 \; F'(u^n)  \; \tau(u^n) 
 = O(\ep \dt^2). \]
Putting this together, we have a one-step error 
\[ E = E_{sch} + E_{per} = O(\dt^3) +  O(\ep \dt^2),\]
so that over the course of the simulation we expect to see error of the form
 \[ Error = O(\dt^2) + O(\epsilon \dt ),\] so that we expect to see   second order results as long as
 $\epsilon \dt $ is small enough, and after that will produce results that look like $\epsilon \dt $.
This explains  the excellent  convergence we observe  in Figure \ref{fig:IMR1}.

\section{Efficient mixed-precision Runge--Kutta methods} \label{sec:methods}
In this section we exploit the framework in Section \ref{sec:framework}  to develop a mixed precision approach to Runge--Kutta methods.
We first show how we can add correction steps into a naive implementation of a mixed-precision methods to raise the
perturbation order of the method, as computed by the conditions in Section \ref{sec:conditions}.
Next we use the order and perturbation conditions to develop novel efficient methods that have high  consistency order and 
high perturbation order using an appropriate optimization code similar to those described in \cite{KetchRKOpt}.

\subsection{Mixed precision implementation and corrections to known Runge--Kutta methods}
In this section, we show that often,  low-precision computation of the  implicit function yields naive mixed-precision methods
that have perturbation errors that may degrade the accuracy of the solution for sufficiently small $\dt$.
It is possible to correct this by adding high order explicit steps; this approach yields methods that can be shown to
satisfy both the consistency \eqref{OC} and perturbation conditions \eqref{PC}.

\subsubsection{Implicit Midpoint rule with correction}
The mixed precision implicit midpoint rule we described above \eqref{IMPButcher_MP}
\begin{align*} 
y^{(1)}&= u^n + \frac{\Delta t}{2} \Fep\left(y^{(1)} \right)\\
u^{n+1}&= u^n + \Delta t F \left(y^{(1)} \right) .
\end{align*}
has global error
\[ Error = O(\dt^2) + O( \epsilon \dt) \]
so that it gives  second order ($O(\dt^2)$)  results as long as $\epsilon \dt $ is small enough; once $\dt$ gets small compared to $\ep$ we observe
degraded convergence.  To eliminate this error, we wish to modify the method   \eqref{IMPButcher_MP} 
so that the $O(\ep \dt^2)$ term in the expansion \eqref{expansion} is set to zero. The framework above suggests how this can be done.
To improve the order of convergence, when $\dt$ is small we add correction terms into the  mixed-precision method:
\begin{subequations} \label{IMPButcher_MPcorr} 
\begin{align} 
y^{(1)}_{[0]} &= u^n+ \frac{1}{2} \dt  \Fep(y^{(1)}_{[0]}) \label{IMPButcher_MPcorr1} \\ 
y^{(1)}_{[k]} & = u^n+ \frac{1}{2}  \dt  F(y^{(1)}_{[k-1]})  \; \; \; \mbox{for} \; \; k = 1,...,  p-1  \label{IMPButcher_MPcorr2}\\
u^{n+1}&= u^n + \Delta t F \left(y^{(1)}_{[p-1]})  \right) . \label{IMPButcher_MPcorr3}
\end{align}
\end{subequations}
so that
\[ A = \left( \begin{array}{ll}
0 & 0  \\
\frac{1}{2} & 0   \\
\end{array} \right), \; \; \; 
c=  \left( \begin{array}{l}
0 \\ \frac{1}{2} \\ 
\end{array} \right), \; \; \; 
b=  \left( \begin{array}{l}
0 \\  1 \\
 \end{array} \right), \; \; \; 
\]
\[ \Aep = \left( \begin{array}{lll}
\frac{1}{2} & 0  \\
0 & 0 \\
\end{array} \right), \; \; \; 
\cep=  \left( \begin{array}{l}
\frac{1}{2} \\ 0\\
 \end{array} \right), \; \; \; 
\bep=  \left( \begin{array}{l}
0 \\ 0 \\
 \end{array} \right), \; \; \; 
\]
\[ \At =A+ \Aep = \left( \begin{array}{lll}
\frac{1}{2} & 0   \\
\frac{1}{2} & 0  \\
\end{array} \right), \; \; \; 
\ct=c+ \cep=  \left( \begin{array}{l}
 \frac{1}{2}  \\ \frac{1}{2}  \\
 \end{array} \right), \; \; \; 
\bt=b+\bep=  \left( \begin{array}{l}
 0 \\ 1\\
  \end{array} \right). \; \; \; 
\]
We observe that to zero out the  $O(\ep \dt^2)$ term we require
\[ \bep \ct=0 , \; \;  \bt \cep =  0 .\]
The method \eqref{IMPButcher_MPcorr} satisfies these equations, and so we obtain global error:
\[ Error = O(\dt^2) + O(\ep \dt^2).\]

Note that this approach to reduce precision errors is reminiscent of that in \cite{Kouya2012}; the framework we developed 
allows us to understand this correction approach as a new method.

\begin{figure}[H]
     \begin{center}
\includegraphics[width=0.325\textwidth]{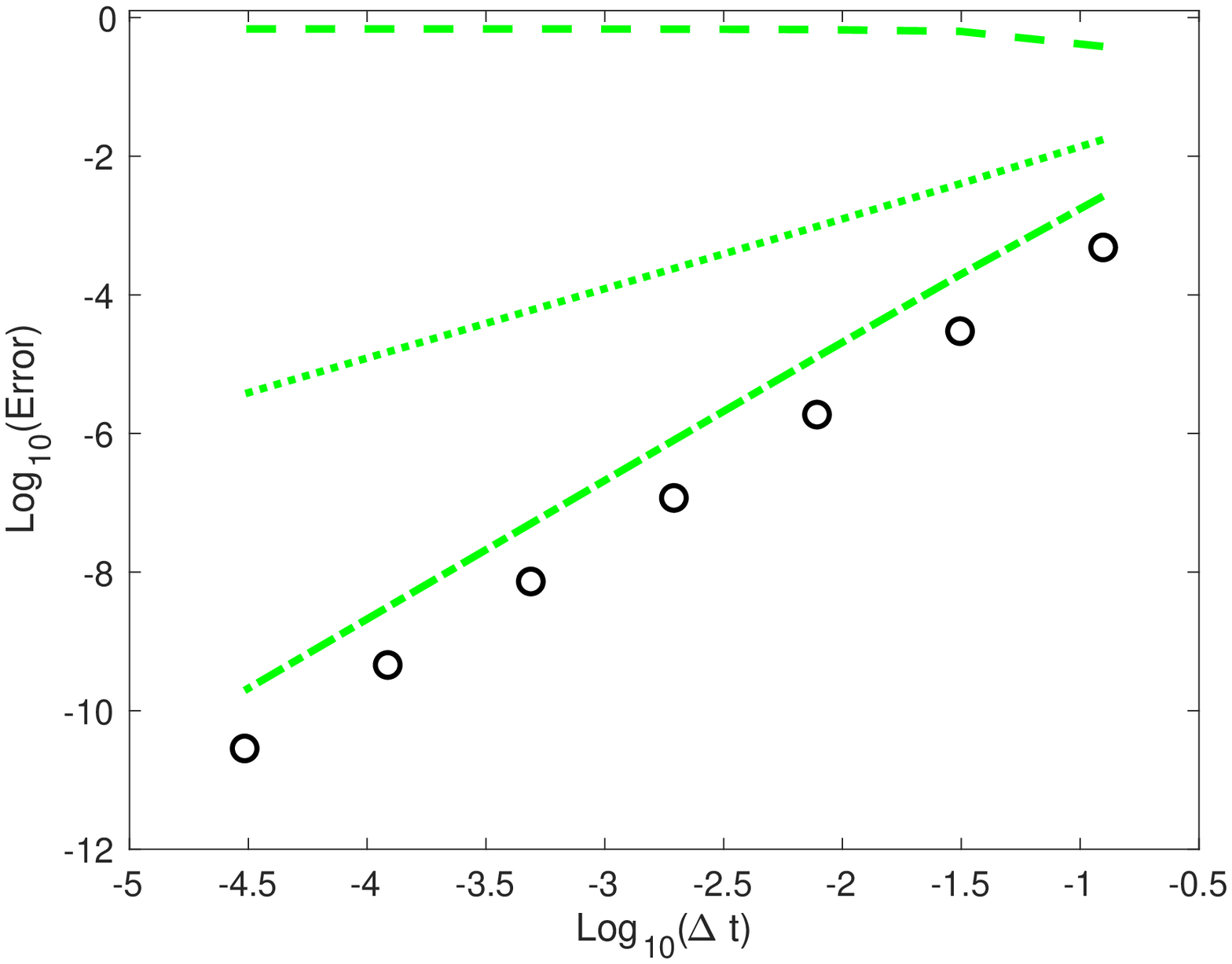}
\includegraphics[width=0.325\textwidth]{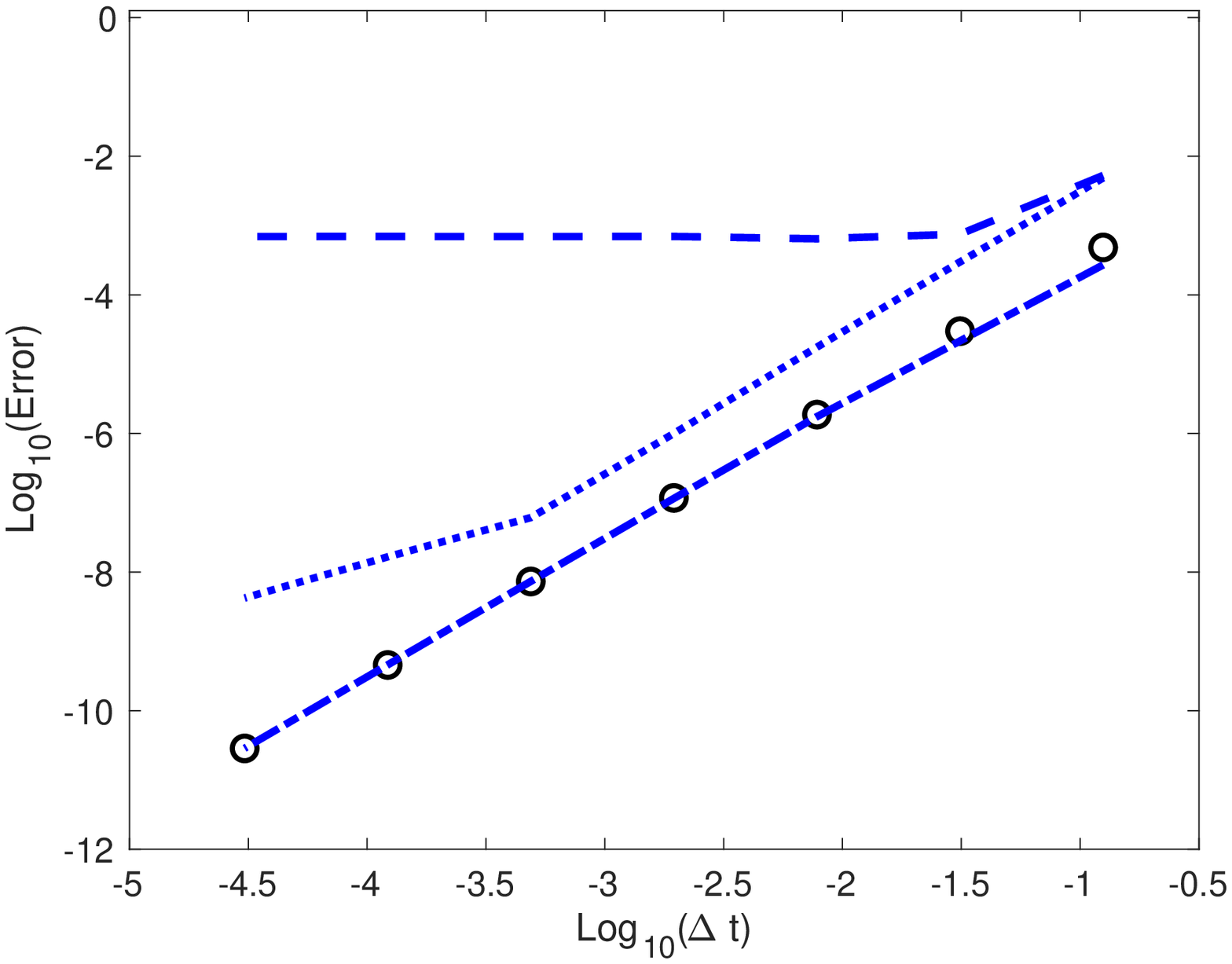}
\includegraphics[width=0.325 \textwidth]{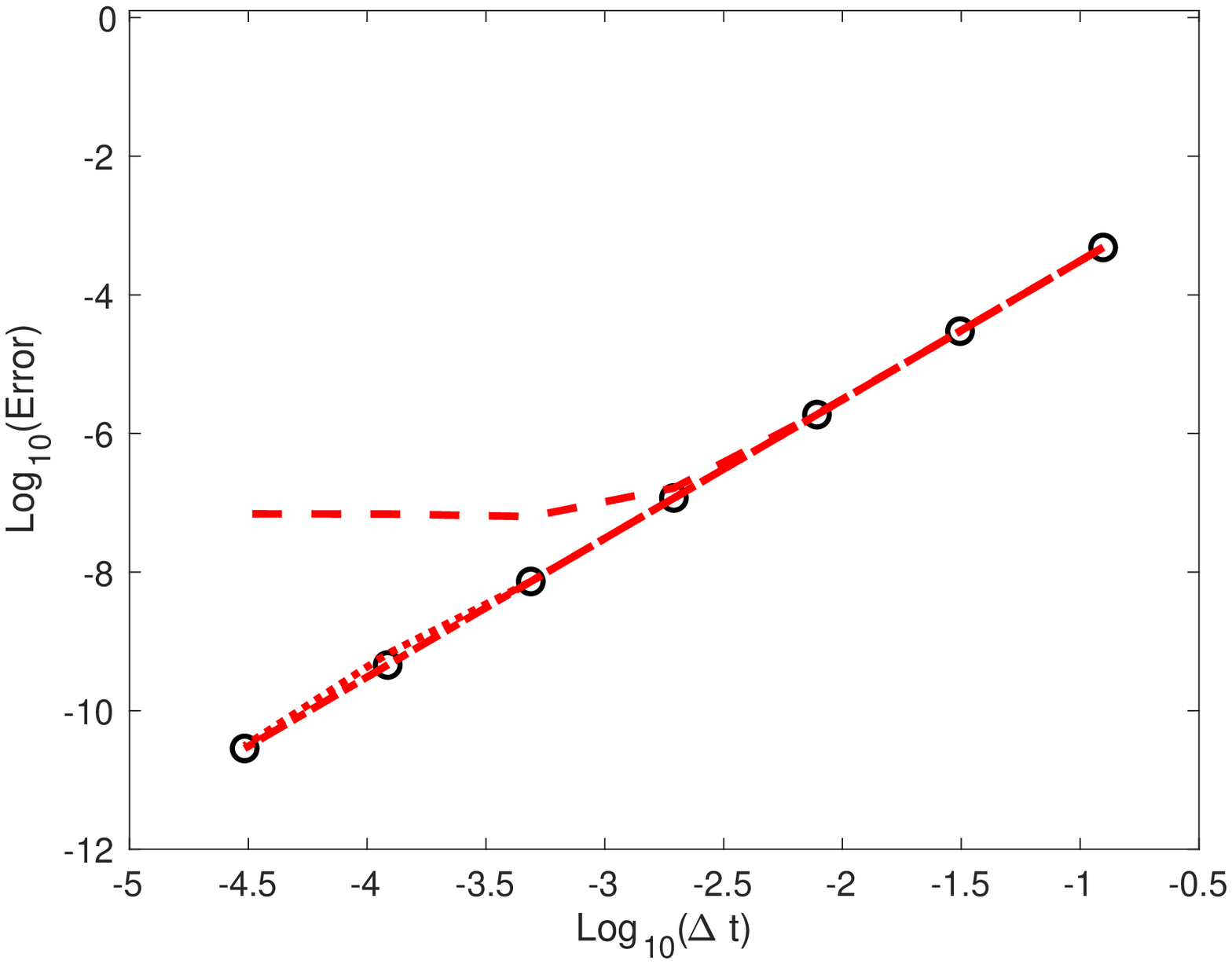}\\
\caption{The implicit midpoint rule applied to the van der Pol system \eqref{eqn:vdp}
with correction steps. Left: Zero precision $\ep=O(1)$, Middle: half precision  $\ep=O(10^{-4})$, Right: single precision   $\ep=O(10^{-8})$.
}
\label{fig:IMRsweep}
   \end{center}
\end{figure}

\noindent{\bf Numerical Results:} 
In the following we demonstrate how this  method performs in practice. 
As before, we use the  van der Pol system, Equation \eqref{eqn:vdp}
  with $a=1$ and initial conditions $y_1(0)= 2$ and $y_2(0) = 0$. We stepped this forward 
  using the implicit midpoint rule  \eqref{IMPButcher_MPcorr} to a final time $T_f=1.0$.
  We show how  using a mixed precision implementation and then a correction step \eqref{IMPButcher_MPcorr2} improves the error.
  In Figure \ref{fig:IMRsweep} we show the results for  $\ep=O(1)$ (zero precision, left), $\ep=O(10^{-4})$ (half precision), and $\ep=O(10^{-8})$ (single precision, right).
  The low precision  \eqref{IMPButcher_LP} is shown in a dashed line, the mixed precision method  \eqref{IMPButcher_MPcorr} with no correction ($p=1$)
 in a dotted line,  and the mixed precision method  \eqref{IMPButcher_MPcorr} with one  correction $(p=2)$ in a dash-dot line. 
 The reference solution  computed in double precision is shown in black circle markers.
   
   Figure  \ref{fig:IMRsweep} shows that each progressive correction, which involves explicit computation of the function, produces a more accurate method. 
   Looking at the perturbation conditions, we showed in Section \ref{IMRexplained}  that the  low precision method  \eqref{IMPButcher_LP} has errors  
   $ Error_{LP} =  O(\dt^2) + O( \epsilon).$  As expected, dashed line solutions shown in Figure \ref{fig:IMRsweep} all have flat line errors at the level 
   of their respective $\epsilon$ values.
   
   The mixed precision  method  \eqref{IMPButcher_MPcorr} with $p=1$ (this is the same method give by Equation \eqref{IMPButcher_MP}), has errors
     $  Error_{p=1} =  O(\dt^2) + O( \epsilon \dt),$  which is reflected in the fact that the dotted line solutions  for the zero precision case in Figure \ref{fig:IMRsweep} 
      has slope  $\sigma =1$. For the half precision case the dotted line solution starts off with a slope $\sigma =2$, but once the time-step gets sufficiently small
      we see the line changes and now has slope $\sigma =1$: this shows clearly that once $\dt$ is small enough compared to $\epsilon$, the $O( \epsilon \dt)$ 
      term dominates and we see first order convergence. For the single precision, we don't observe this phenomenon in this example because $\dt$ is not small enough 
      compared to $\epsilon$.
   
   Finally, the mixed precision method with one correction step (Equation  \eqref{IMPButcher_MPcorr} with $p=2$ ) has errors
     $ Error_{p=2} =  O(\dt^2) + O( \epsilon \dt^2).$
   For the half and single precision,  two corrections steps produce a second order solution, and this line has $\sigma = 2$. 
 The order  of the same accuracy as a complete computation in double precision (shown in black circles). 
This would be worth-while in all cases where  two explicit steps take much less computational time than the savings realized from 
replacing a double-precision implicit solve with one that has single or half precision. The case of zero precision, $\ep=O(1)$,  
has a larger error, but the correction step  clearly gives an error with slope $\sigma=2$.

\subsubsection{A mixed precision implementation of a two stage  third order SDIRK}

We take the two stage  third order singly diagonally implicit method \cite{Alexander1977} 
\begin{subequations}  \label{SDIRK}
\begin{eqnarray}     
y^{(1)} &=& u^n + \gamma \dt F(y^{(1)})  \\
     y^{(2)} &=& u^n + \left( 1- 2 \gamma \right)  \dt F(y^{(1)})  + \gamma \dt F(y^{(2)})   \\
      u^{n+1}&=&u^n+ \frac{1}{2} \dt F(y^{(1)})  + \frac{1}{2} \dt F(y^{(2)}) 
\end{eqnarray}
\end{subequations}
with $\gamma =  \frac{\sqrt{3}+3}{6}$.
A low-precision implementation is given by
\begin{subequations} \label{SDIRK-LP}
\begin{eqnarray}
     y^{(1)} &=& u^n + \gamma \dt \Fep(y^{(1)})  \\
     y^{(2)} &=& u^n + \left( 1- 2 \gamma \right)  \dt \Fep(y^{(1)})  + \gamma \dt \Fep(y^{(2)})   \\
      u^{n+1}&=&u^n+ \frac{1}{2} \dt F(y^{(1)})  + \frac{1}{2} \dt \Fep(y^{(2)}) .
\end{eqnarray}
\end{subequations}
Using a low-precision computation only of the  implicit function yields the mixed-precision method:
\begin{subequations} \label{SDIRK-MP}
\begin{eqnarray}
     y^{(1)} &=& u^n + \gamma \dt \Fep(y^{(1)})  \\
     y^{(2)} &=& u^n + \left( 1- 2 \gamma \right)  \dt F(y^{(1)})  + \gamma \dt \Fep(y^{(2)})   \\
      u^{n+1}&=&u^n+ \frac{1}{2} \dt F(y^{(1)})  + \frac{1}{2} \dt F(y^{(2)}) .
\end{eqnarray}
\end{subequations}
The consistency conditions are satisfied to order three. The highest order non-zero perturbation term is
$ \bt \cep = \gamma $
 so we have a perturbation error $E_{per} = O(\dt^2 \ep)$ at each time step, or a global error
 \[ Error = O(\dt^3) + O(\ep \dt  ) .\]
 When $\dt^2 < \ep$, the perturbation error will dominate. In this case, we can correct the method
 by adding explicit stages:
 \begin{subequations} \label{SDIRK-MPfix}
\begin{eqnarray}
y_{[0]}^{(1)} &=& u^n+\gamma \dt \Fep(y^{(1)}_{[0]}) \\
y_{[k]}^{(1)}  &=& u^n+\gamma \dt  F(y^{(1)}_{[k-1]}) \; \; \; \mbox{for} \; \; k = 1,...,  p-1 \\
y_{[0]}^{(2)}  &=& u^n+(1-2\gamma) \dt F(y^{(1)}_{[p-1]})+\gamma \dt F(y^{(2)}_{[0]}) \\
y_{[k]}^{(2)}  &=& u^n+(1-2\gamma) \dt F(y^{(1)}_{[p-1]})+\gamma \dt F(y^{(2)}_{[k-1]}) \; \; \; \mbox{for} \; \; k = 1,...,  p-1 \\
  u^{n+1}&=&u^n+ \frac{1}{2} \dt F(y^{(1)}_{[p-1]}))  + \frac{1}{2} \dt F(y^{(2)}_{[p-1]})) .
  \end{eqnarray}
\end{subequations}
This corrected method with $p=3$ has coefficients:
\[ A= \left( \begin{array}{cccccc}
0 & 0 & 0 & 0 & 0 & 0 \\
\gamma & 0 & 0 & 0 & 0 & 0 \\
0 & \gamma & 0 & 0 & 0 & 0  \\
0 & 0 & (1-2\gamma)  &0 & 0 & 0 \\
0 & 0 & (1-2\gamma)  & \gamma & 0 & 0 \\
0 & 0  & (1-2\gamma)  & 0 &  \gamma & 0  \\
\end{array} \right), \; \; \; \;
\Aep= \left( \begin{array}{cccccc}
\gamma & 0 & 0 & 0 & 0 & 0  \\
0 & 0 & 0 & 0 & 0 & 0 \\
0 & 0 & 0 & 0 & 0 & 0 \\
0 & 0 & 0 & \gamma & 0 & 0 \\
0 & 0 & 0 & 0 & 0 & 0 \\
0 & 0 & 0 & 0 & 0 & 0 \\
\end{array} \right) \]
\[
b = \left( \begin{array}{cccccc}
0 & 0 & \frac{1}{2} & 0 & 0 & \frac{1}{2} \\
\end{array} \right) ,\; \; \; \;
\bep = \left( \begin{array}{cccccc}
0 & 0 & 0 & 0 & 0 & 0 \\
\end{array} \right) .\;
\]
The order conditions are, as before, satisfied to third order, and the perturbation terms contribute to the errors terms of the form $O(\ep \dt^4)$, 
so the overall global error from the method with three corrections for each implicit solve is 
\[Error_{p=3}  = O(\dt^3) +  O(\ep \dt^3).\]

\noindent{\bf Numerical Results:} 
To demonstrate the performance of this  method in practice, we apply its various implementations to 
the  van der Pol system, Equation \eqref{eqn:vdp}, with $a=1$ and initial conditions $y_1(0)= 2$ and $y_2(0) = 0$. 
We step this forward to a final time $T_f=1.0$,  using low precision, mixed precision, and mixed precision with successive corrections.

    \begin{wrapfigure}[21]{r}{0.5\linewidth} \vspace{-0.25in}
    \begin{center}
\includegraphics[width=0.485\textwidth]{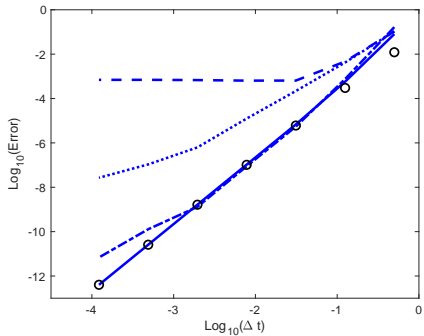}
\caption{The SDIRK method  with successive corrections  applied to the van der Pol system \eqref{eqn:vdp}.
The half-precision implementation \eqref{SDIRK-LP})  shown as a dashed line. The mixed precision implementation
\eqref{SDIRK-MP}) with $\ep= O(10^{-4})$ is shown with a dotted line. The mixed precision implementation \eqref{SDIRK-MPfix} with  $\ep= O(10^{-4})$
and one ($p=2$) correction steps is shown with dash-dot, and with two ($p=3$) correction steps with a  solid line. 
The double precision implementation of \eqref{SDIRK} is shown for reference with black dots.
}
\label{fig:SDIRKsweep}
\end{center}
\end{wrapfigure} 
In Figure \ref{fig:SDIRKsweep} we show the half-precision results of the various implementations of the SDIRK method.
 First, we use a half-precision implementation of  the SDIRK method, as given in Equation  \eqref{SDIRK-LP}) with $\ep= O(10^{-4})$.
 The errors resulting from this implementation are shown by a dashed line. This line is horizontal at the level of $O(\ep)$, as expected from 
  the error given by this implementation: $Error_{LP} =  O(\dt^3) + O( \epsilon). $
Using the naive mixed precision  implementation \eqref{SDIRK-MPfix} we expect an error of 
$Error_{p=1} =  O(\dt^3) + O( \epsilon \dt). $ In  Figure \ref{fig:SDIRKsweep}, the dotted line shows that error initially has slope $\sigma=2$:
this happens when $\ep$ is small compared to  $\dt$ and so the $O( \epsilon \dt)$ looks like  $O( \dt^2)$. 
However, as $\dt$ gets smaller, we see the error line become first order.  Adding one correction step to the mixed precision implementation 
\eqref{SDIRK-MPfix}, we obtain an error shown as the  dash-dot line, which initially looks third order (slope $\sigma=3$)
but then, as $\dt$ becomes small compared to $\epsilon$, begins to look like it is second order (slope $\sigma=3$). This matches the expected order
$Error_{p=2} =  O(\dt^3) + O( \epsilon \dt^2). $ Finally, when we add two correction steps to the mixed precision implementation, the error (shown 
as a solid line) has slope $\sigma=3$, as expected from the predicted order 
$Error_{p=3} =  O(\dt^3) + O( \epsilon \dt^3). $  This solution matches the double precision reference solution  shown in black circle markers.

\subsubsection{Mixed precision implementation of a two stage L-stable scheme}
Consider the L-stable fully implicit Lobatto IIIC scheme \cite{NorsettWanner1981}:
\begin{subequations} \label{Lobatto}
 \begin{eqnarray}
  y^{(1)} &=& u^n +\frac{1}{2} \Delta t F(  y^{(1)} )- \frac{1}{2} \Delta t  F( y^{(2)} )\\
  y^{(2)}  &=& u^n +\frac{1}{2} \Delta t F(  y^{(1)} ) + \frac{1}{2} \Delta t  F( y^{(2)} )\\
    u^{n+1}  &=& y^{(2)}.
  \end{eqnarray}
  \end{subequations} 
A naive mixed precision implementation of this method is given by:
\begin{subequations} \label{Lobatto-MP}
 \begin{eqnarray}
  y^{(1)} &=& u^n +\frac{1}{2} \Delta t \Fep(  y^{(1)} )- \frac{1}{2} \Delta t  \Fep( y^{(2)} )\\
  y^{(2)}  &=& u^n +\frac{1}{2} \Delta t \Fep(  y^{(1)} ) + \frac{1}{2} \Delta t  \Fep( y^{(2)} )\\
    u^{n+1}  &=& u^n +\frac{1}{2} \Delta t F(  y^{(1)} ) + \frac{1}{2} \Delta t  F( y^{(2)} ) 
  \end{eqnarray}
  \end{subequations} 
  This implementation has Butcher coefficients
  \[ \Aep = \left( \begin{array}{rr}
\frac{1}{2}  & - \frac{1}{2}   \\
\frac{1}{2}  & \frac{1}{2} \\
\end{array} \right),  \; \; \; \; 
A= \left( \begin{array}{rr}
0 & 0    \\
0 & 0 \\
\end{array} \right), \]
\[ \bep= \left( \begin{array}{ll} 0 & 0  \end{array} \right), \; \; \; \; 
b=\left( \begin{array}{lll}  \frac{1}{2} & \frac{1}{2} \end{array} \right).\]
The error coming from this implementation is, according to the analysis in Section \ref{sec:framework}, 
\[ Error = O(\dt^2) + O(\ep \dt).  \]
  
A corrected mixed precision implementation is given by
\begin{subequations} \label{Lobatto-MPfixed}
 \begin{eqnarray}
  y^{(1)}_{[0]} &=& u^n +\frac{1}{2} \Delta t \Fep(  y^{(1)} )- \frac{1}{2} \Delta t  \Fep( y^{(2)} )\\
  y^{(2)}_{[0]}  &=& u^n +\frac{1}{2} \Delta t \Fep(  y^{(1)} ) + \frac{1}{2} \Delta t  \Fep( y^{(2)} )\\
y^{(1)}_{[1]} &=& u^n +\frac{1}{2} \Delta t F(  y^{(1)}_{[0]} )- \frac{1}{2} \Delta t  F( y^{(2)}_{[0]} )\\
  y^{(2)}_{[1]}  &=& u^n +\frac{1}{2} \Delta t F(  y^{(1)}_{[0]}) + \frac{1}{2} \Delta t  F( y^{(2)}_{[0]} )\\
    u^{n+1}  &=& u^n +\frac{1}{2} \Delta t F(  y^{(1)}_{[1]} ) + \frac{1}{2} \Delta t  F( y^{(2)}_{[1]} ) 
    \end{eqnarray}
  \end{subequations}

\[ \Aep = \left( \begin{array}{rrrr}
\frac{1}{2}  & - \frac{1}{2}  &  0 & 0 \\
\frac{1}{2}  & \frac{1}{2} &  0 & 0 \\
0  & 0 &  0 & 0 \\
0 & 0 &  0 & 0 \\
\end{array} \right),  \; \; \; \; 
A= \left( \begin{array}{rrrr}
0 & 0 &  0 & 0 \\
0 & 0 &  0 & 0 \\
\frac{1}{2}  & - \frac{1}{2}  &  0 & 0 \\
\frac{1}{2}  & \frac{1}{2} &  0 & 0 \\
\end{array} \right), \]
\[ \bep= \left( \begin{array}{llll} 0 & 0 & 0 & 0 \end{array} \right), \; \; \; \; 
b=\left( \begin{array}{lllll} 0 & 0 & \frac{1}{2} & \frac{1}{2} \end{array} \right).\]
Using the analysis in Section \ref{sec:framework} we see that the error coming from this implementation is expected to be
\[ Error = O(\dt^2) + O(\ep \dt^3) .\]

\noindent{\bf Numerical Results:} 
As before, we apply the  different possible  implementations  to evolve  the  van der Pol system, Equation \eqref{eqn:vdp},
 (with $a=1$ and initial conditions $y_1(0)= 2$ and $y_2(0) = 0$)  to a final time $T_f=1.0$.
 
    \begin{figure}[H]
     \begin{center}
\includegraphics[width=0.495\textwidth]{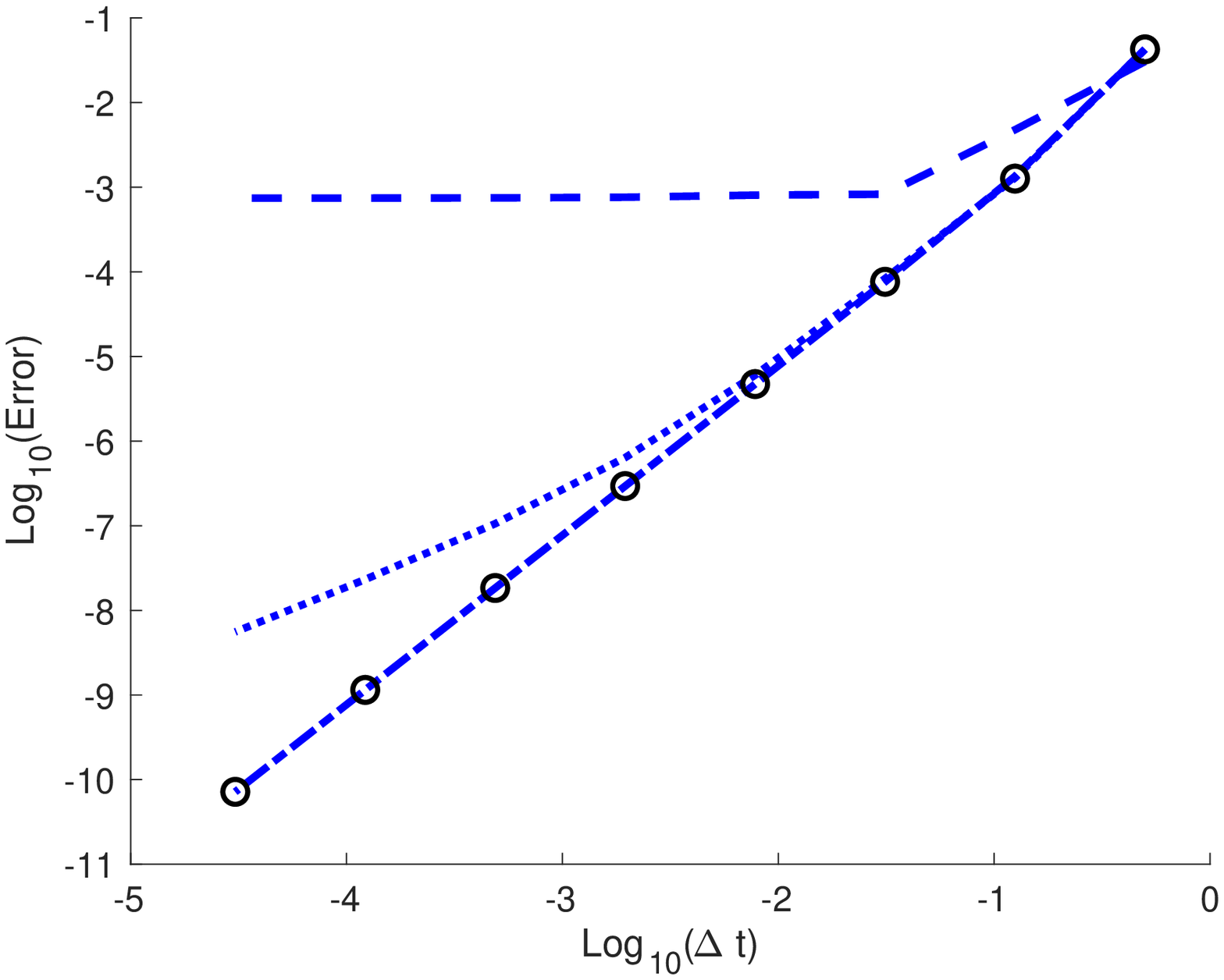}
\includegraphics[width=0.495\textwidth]{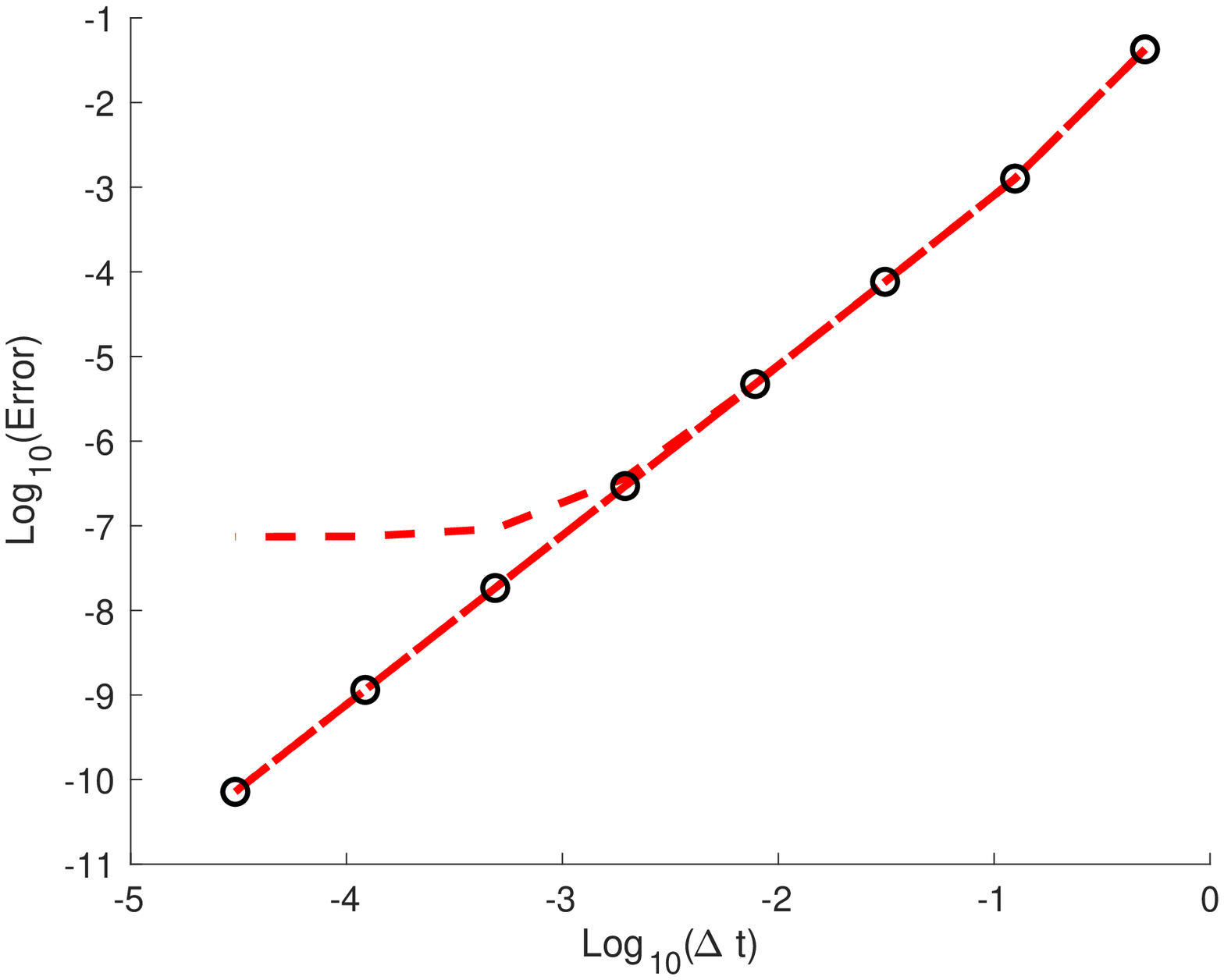}
\caption{The Lobatto method with various implementations. The dashed line is a low-precision 
implementation, the dotted line is the naive mixed precision implementation \eqref{Lobatto-MP},
and the dash-dot line is the corrected mixed precision implementation \eqref{Lobatto-MPfixed}.
Left: half precision  $\ep=O(10^{-4})$, Right: single precision   $\ep=O(10^{-8})$.
}
\label{fig:LobattoFix}
   \end{center}
\end{figure} 

In Figure \ref{fig:LobattoFix} we show the half-precision  (left) and single-precision  (right) results of the various implementations of the Lobatto IIIC method.
The errors from the low precision implementation are shown  by a dashed line which starts off with slope $\sigma =2$ and very quickly 
becomes horizontal, as expected from the error analysis which predicts $Error_{LP} =  O(\dt^2) + O( \epsilon). $
The errors from the naive mixed-precision  implementation \eqref{Lobatto-MP} are shown by a dotted line which has 
slope $\sigma =2$ for   $\ep=O(10^{-8})$, but has slope $\sigma =1$ for  $\ep=O(10^{-4})$, which matches the predicted 
error  $Error_{MP} =  O(\dt^2) + O( \epsilon \dt). $ What is happening here is that for the half precision implementation the term 
$O(\epsilon \dt)$ dominates early on, whereas for the single precision implementation we observe the $O(\dt^2) $ convergence because 
$\dt$ is large compared to $\ep$. (For the single precision implementation the dotted line is hidden by the dash-dot line).
 Finally, the corrected mixed-precision method  \eqref{Lobatto-MPfixed}
 has errors (shown in a dash-dot line) that are second order and that match the full-precision implementation (black circles).

\subsection{Third order novel methods}
\begin{figure}[H]
    \begin{center}
\includegraphics[width=0.485\textwidth]{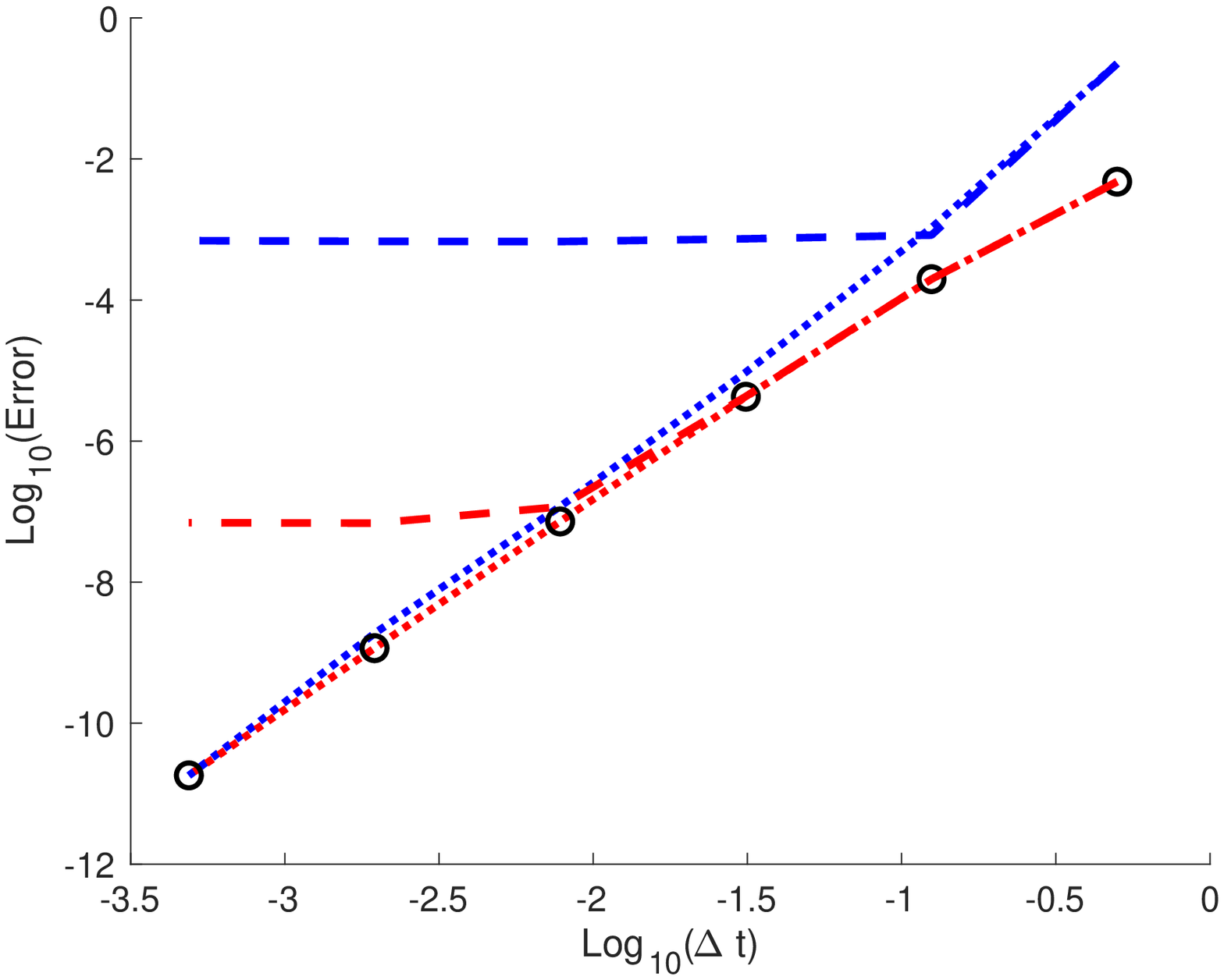}
\includegraphics[width=0.485\textwidth]{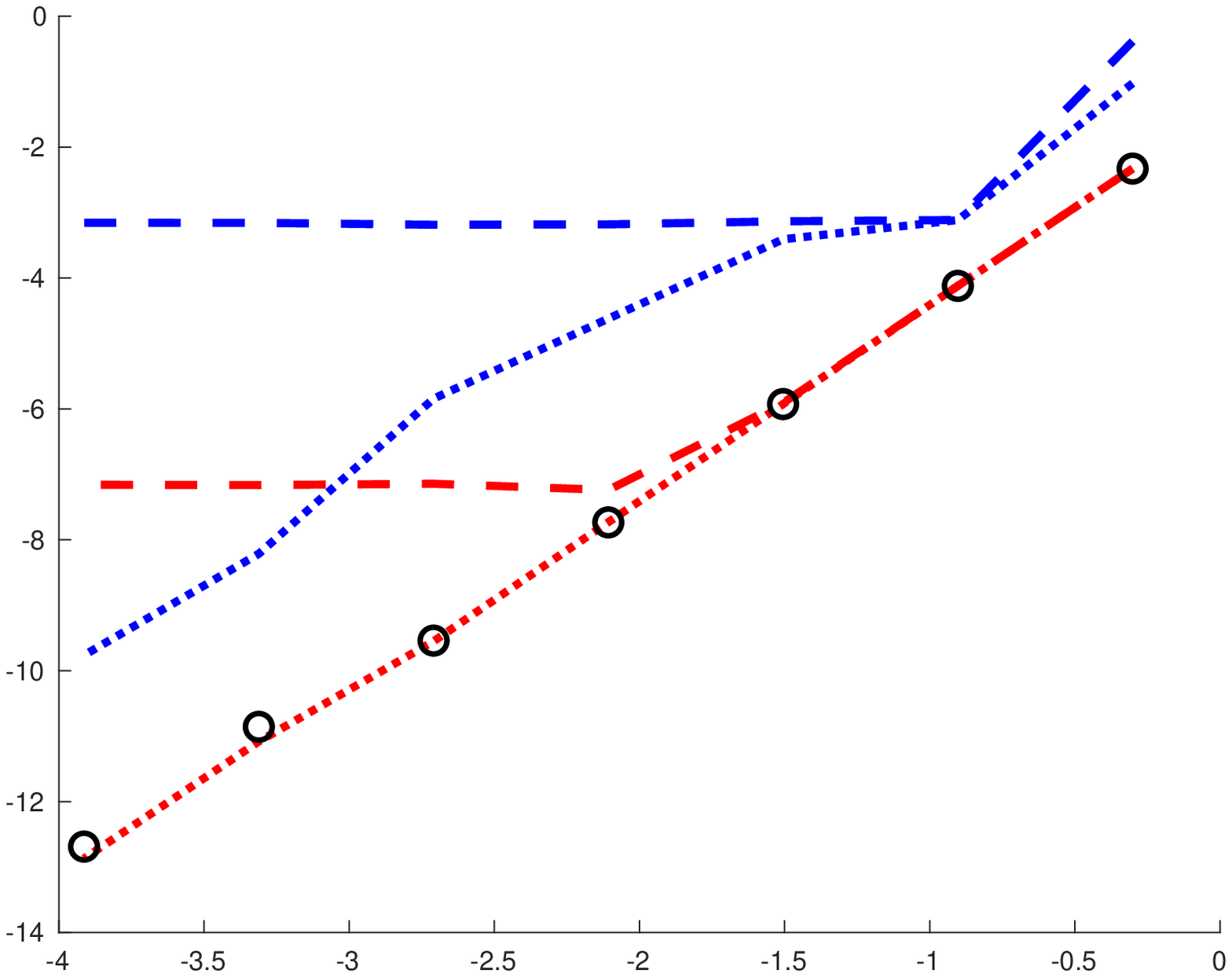}
\caption{The four stage third order method given in Sections \ref{4s3pMPa} (Method 4s3pA, left) and 
 \ref{4s3pMPb} (Method 4s3pB, right).
The half-precision implementation ($\ep= O(10^{-4})$) is shown as a blue dashed line 
and the single-precision implementation  $\ep= O(10^{-8})$ as a red dashed line. 
The mixed precision implementation with $\ep= O(10^{-4})$ is shown with a dotted blue line, 
and the mixed precision implementation with $\ep= O(10^{-8})$ is shown with a dotted red line. 
The double precision implementation is shown for reference with black dots.
}
\label{fig:4s3pMPComparison}
\end{center}
\end{figure}

The framework presented in Section \ref{sec:framework} can not only be used to analyze naive mixed-precision implementations of
existing methods and corrections to such methods, but to devise new methods. In this section we present examples of two methods
that were developed to satisfy the order \eqref{OC} and perturbation conditions \eqref{PC} in Section \ref{sec:conditions} to high order. Both methods are 
four-stage third order methods. The first method, presented in Subsection \ref{4s3pMPa}, is not A-stable, and is third order with high order perturbation errors:
\[ Error = O(\dt^3) + O(\ep \dt^3).\]
The second method, presented in Subsection \ref{4s3pMPb}, is a perturbation of the four-stage third order L-stable method  in \cite{HairerWanner1991},
and so is A-stable. However, its perturbation errors are not as high order:
\[ Error = O(\dt^3) + O(\ep \dt^2).\]
The difference between these methods is evident in Figure \ref{fig:4s3pMPComparison}. The errors for the mixed precision implementation of
Method 4s3pA (Figure \ref{fig:4s3pMPComparison}, left) are shown in dotted lines (blue for for half precision, and red for single precision).
The mixed precision errors match the corresponding low-precision errors initially, but as $\dt$ gets smaller, the mixed
precision errors match with the double-precision errors. For the A-stable Method 4s3pB (Figure \ref{fig:4s3pMPComparison}, right)
 this is also true when considering the mixed  precision method with $\ep= O(10^{-8})$. However, the mixed  precision method with $\ep= O(10^{-4})$
 does not match the double precision errors, even as $\dt$ gets small.


\subsubsection{A four-stage third order mixed-precision method} \label{4s3pMPa}
This method, referred to as {\bf Method 4s3pA} is given by the following coefficients:

The matrix $A$ is
\[ a_{21}=0.211324865405187, \; \;
a_{31}=0.709495523817170, \; \;
a_{32}=-0.865314250619423,\]
\[ a_{41}= 0.705123240545107, \; \;
a_{42} =0.943370088535775, \; \;
a_{43}= -0.859818194486069,\]
$\Aep$ has coefficients
\[a^\ep_{11}= 0.788675134594813,  \; \;
a^\ep_{31}= 0.051944240459852, \; \;
a^\ep_{33} = 0.788675134594813,\] 
and the vectors are given by
\[ b= ( 0, \frac{1}{2}, 0 , \frac{1}{2}) , \; \; \; \; 
\bep = ( 0, 0, 0 , 0) .\]


\subsubsection{A four-stage third order A-stable mixed-precision method}  \label{4s3pMPb}
This A-stable method, referred to as {\bf Method 4s3pB} is given by the following coefficients:
For this A-stable method, $A$ is given by 
{\small
\[ a_{21}= 2.543016042796356  \; \; 
 a_{31}=   2.451484396921318 , \; \;   a_{32} = 0.024108961241221, \]             
\[ a_{41}=   2.073861819468268   \; \;   a_{42} =  2.367724727682735 \; \;   a_{43} =    1.711868223075524,\]}
 $\Aep$ is 
{\small
 \[ a^\ep_{11}=   a^\ep_{22} = a^\ep_{33} = a^\ep_{44}=    0.5,\]
\[ a^\ep_{21}= -2.376349376129689,  \; \; 
 a^\ep_{31}=   -2.951484396921318 ,  \; \;  a^\ep_{32}=  0.475891038758779,\]
\[  a^\ep_{41}=   -0.573861819468268, \; \;    a^\ep_{42}=  -3.867724727682735, \; \; 
a^\ep_{43}=   -1.211868223075524,\]} 
and  \[ b = \left( \frac{3}{2}, - \frac{3}{2}, \frac{1}{2}, \frac{1}{2} \right), \; \; \; \;  b^\ep = \left( 0, 0 , 0, 0 \right).\]
This method is a perturbation  the four-stage third order L-stable method  in \cite{HairerWanner1991}.

\subsection{A method that satisfies the simplified order conditions}
In all of the above we considered methods that satisfied the perturbation conditions \ref{PC}, that apply 
even when $\tau$ is not a well-behaved function. 
In this section, we devise a method that satisfies the less restrictive order conditions \ref{PCsimplified},
that apply only when $\tau$ is well-behaved. This method,  {\bf Method 4s3pC}  is given by the coefficients
 {\small
\[ a_{21}= -0.050470366527530, \; \;
a_{31}=   0.368613367355336, \; \; a_{32}= 0.273504374252976,\]
\[   a_{41}= 1.803794668975043, \; \;  a_{42}=   0.097485042980759, \; \;  a_{43}= -1.895660952342050.\]
 \[ a^\ep_{11}=  0.511243008730995, \;  \; 
  a^\ep_{21}=   -1.999347282862640, \; \; 
 a^\ep_{22}= 1.957161067302390, \]
 \[      a^\ep_{31}=0.443312893511937, \; \;
   a^\ep_{32}= -0.573131033672219, \; \;
    a^\ep_{33}=   0.128283796414019,\]
     \[  a^\ep_{42}=-0.160330320741428, \; 
     a^\ep_{43}=   0.579597314161362, \; 
     a^\ep_{44}=   1.484688928981990,\]
       \[      a^\ep_{41}= -2, \; \; \; b^\ep = (0, 0, 0, 0),\]
 \[b=(0.002837446974069, \; \;    0.336264433650450, \; \;   0.806376720267787, \; \;   -0.145478600892306).
  \]
}

When we have a well-behaved $\tau$, this method gives errors of the form
\[ Error = O(\dt^3) + O(\ep \dt^3).\]
However, $\tau$ is not well-behaved, this method gives error of the form
\[ Error = O(\dt^3) + O(\ep \dt^2).\]

 \begin{wrapfigure}[26]{r}{0.55\linewidth} \vspace{-0.3in}
\includegraphics[width=0.55\textwidth]{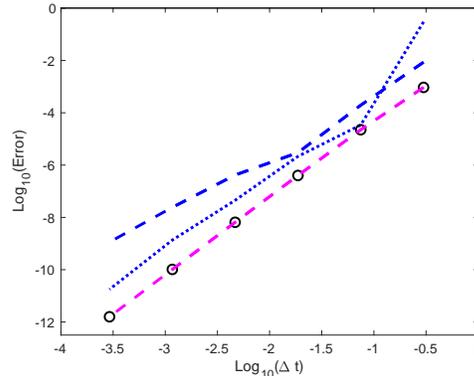}
\caption{{Methods {\bf 4s3pC} and {\bf 4s3pA} for time evolution of  the diffusion equation.
The magenta dashed  line is {\bf Method 4s3pC} with $F$ resulting from 
a Fourier spectral method approximation to $u_{xx}$, and $\Fep$ resulting from 
a second order centered difference approximation to $u_{xx}$.
The dashed blue line is {\bf Method 4s3pC} with $F$ resulting from a double precision 
implementation of a Fourier spectral method approximation to $u_{xx}$,
and $\Fep$ resulting from a half precision "chopping" of the  Fourier spectral method approximation.
The dotted blue line is {\bf Method 4s3pA}  for this same scenario.
The double precision implementation with a Fourier spectral method approximation
 is shown for reference with  black circles.
}}
\label{fig:CompareTau}
\end{wrapfigure} We test this problem on the diffusion equation
\[ u_t = u_{xx} \]
on $x \in (0, 2 \pi)$ with initial conditions $u(x,0)=sin(x) $ and periodic boundaries.

To discretize the spatial derivative we use a high resolution Fourier spectral method for $F$.
For the $\Fep$, we consider two different approaches. In the first approach, we use the low resolution 
centered difference scheme for $u_{xx}$ to evaluate $\Fep$. This is a highly sensitive process,
and a careful stability analysis must be carried out with the two different operators, 
so we do not  recommend trying this approach  in general without rigorous justification. 
We use it here only to  illustrate the effect of using different resolution methods in our
perturbed Runge--Kutta framework. We show in Figure \ref{fig:CompareTau}, 
that using {\bf Method 4s3pC} on the method where $\Fep$ is given by a centered difference
approximation (i.e. a well-behaved $\tau$)  results, as expected, in an error of  
$Error = O(\dt^3) + O(\ep \dt^3)$  (magenta dashed line).

For comparison, we show a low  precision (using the chop command)
approximation of the Fourier spectral method for $\Fep$. 
Recall that  {\bf Method 4s3pC} was designed to work with a well behaved $\tau$. In the
case where we approximate $u_{xx}$ with the Fourier spectral method for $F$ and
with  the centered difference scheme for $\Fep$, we can show that the difference between these is a well-behaved function.
Using {\bf Method 4s3pC} with the mixed precision approach with $\ep=10^{-4}$ (dashed blue line)
results in an error that is initially third order and reduces to second order as $\dt$ gets smaller. 
This matches our expected error $Error = O(\dt^3) + O(\ep \dt^2)$ when  $\tau$ is not well-behaved. 
Compare this with the performance of {\bf Method 4s3pA}, which was designed to work with badly behaved $\tau$;
we see that these errors start at higher than third order, and settle down to third order behavior. 
The convergence of this mixed precision method is of the same order as the low/high resolution method, and of the high 
resolution method, but it has a larger error constant.

\section{Conclusions} \label{sec:conclusions}
In this work we presented a framework for the error analysis of perturbations of Runge--Kutta methods.
In particular, we investigate the case where perturbations  arise from a mixed precision implementation
of Runge--Kutta methods. This is particularly useful for implicit methods, where is implicit evaluation is 
computationally costly.  Using this framework, we investigate mixed precision implementations 
of existing methods and a  correction approach that improves the errors, and devised new methods
that have favorable scheme error and perturbation error properties. Numerical demonstrations
illustrate the performance of these methods as described by the theory.
Although this  mixed precision approach was designed for implicit Runge--Kutta schemes, 
it can also be applied when repeated  explicit function evaluation is expensive, 
and storage of the computed values  is not possible due to the size of the problem.

In the case where we use a chopping routine to emulate a low precision operator, 
we developed more stringent conditions on the method to handle the unbounded behavior of the truncation operator.
The framework developed holds for more general perturbations than mixed precision calculations: 
we also presented simplified order conditions that are applicable when the perturbation function
$\tau$ is well-behaved. These methods can thus be extended to  many types of perturbations.
While in this work we treat epsilon as a single constant upper bound, in future work we will 
generalize this approach  to design methods with varying orders of epsilon.

\end{document}